\documentclass[a4paper,english]{article}

\usepackage[english]{babel}
\usepackage[utf8]{inputenc}
\usepackage{hyperref}
\usepackage{amsmath}
\usepackage{amsfonts}
\usepackage{amssymb}
\usepackage{mathtools}
\usepackage{amsthm}
\usepackage{graphicx}
\usepackage{tikz}
\usepackage{times}
\usepackage{graphicx}
\usepackage{array}
\usepackage{amsmath}
\usepackage{amsfonts}
\usepackage{amssymb}
\usepackage{mathtools}
\usepackage{enumitem}
\usepackage{todo}
\usepackage{gensymb}
\usepackage{xcolor}
\usepackage{booktabs}
\usepackage{colortbl}
\usepackage{makecell}
\usepackage{hyperref}
\usepackage{comment}
\usepackage[dvipsnames]{xcolor}

\usepackage{booktabs}
\usepackage{siunitx}
\usepackage[margin=1in]{geometry}

\usetikzlibrary{patterns, shapes, arrows, positioning, fit, backgrounds, calc}

\usepackage{zibtitlepage}
\usepackage{authblk}

\definecolor{myblue}{RGB}{7, 127, 247}
\definecolor{myred}{RGB}{203, 65, 84}
\definecolor{mygreen}{RGB}{85, 107, 47}

\newif\ifrevision
\revisiontrue
\DeclareRobustCommand{\rev}[1]{\ifrevision{\color{blue}#1}\else#1\fi}

\revisionfalse

\begin{document}

\title{Risk-Aware Unit Commitment in District Heating Networks: Chance-Constrained and CVaR Optimization Under Demand Uncertainty}

\author[1]{Annika Buchholz\footnote{corresponding author, buchholz@zib.de, \ZTPOrcid{0009-0008-3235-0896}}}
\author[1]{Janina Zittel\footnote{zittel@zib.de, \ZTPOrcid{0000-0002-0731-0314}}}
\author[1,2]{Thorsten Koch\footnote{koch@zib.de, \ZTPOrcid{0000-0002-1967-0077}}
}

\affil[1]{Zuse Institute Berlin, Applied Optimization, Berlin, Germany}
\affil[2]{Technische Universtität Berlin, Software and Algorithms for Discrete Optimization, Berlin, Germany}

\maketitle



\begin{abstract}
While district heating networks are a key component of the energy transition, their operational planning is challenging due to substantial uncertainty in heat demand. We study how unit-commitment schedules can account for this uncertainty and what reliability costs. By adapting chance-constrained programming (CC) and conditional value-at-risk (CVaR) optimization to the mixed-integer unit-commitment problem for district-heating networks, we can compute optimal unit-commitment schedules under demand uncertainty. We generate heat demand scenarios using a Bayesian model, which explicitly quantifies forecasting uncertainty and yields predictive distributions as input to the network optimization. To handle the inherent uncertainty in heat demand, we apply two stochastic optimization approaches that exploit this distribution: chance-constrained programming limits the probability of unmet demand, and in our scenario-based setting reduces to a deterministic MIP in which demand is replaced by an empirical quantile, while CVaR optimization penalizes severe shortfalls in the tail of the distribution. We evaluate the approaches both on real-world data from the Berlin district heating network and on a benchmark set of synthetic, realistically parameterized instances of varying sizes. The results show that the two methods are complementary: CC cheaply traces the high-reliability end of the cost--reliability trade-off and reveals the physical reliability limit of the network, whereas CVaR covers the low-cost, risk-tolerant regime at considerably higher computational effort.
\end{abstract}


\section{Introduction}
\rev{District heating (DH) networks distribute hot water from central generation
units, such as combined heat and power (CHP) plants, boilers, and heat pumps, to
buildings.} Operational planning of \rev{DH} networks centers on
unit commitment (UC): \rev{deciding, for each time step of a planning horizon,
which units are switched on and how much each of them produces. Start-up costs,
minimum up/down times, and ramping limits make these on/off decisions discrete
and coupled over time, so UC is} a mixed-integer \rev{program (MIP)}. The dominant difficulty is
uncertainty in heat demand, driven by weather and consumer behavior. Committing units
against a single point forecast risks costly over-provisioning or, worse, unmet demand. \rev{Instead, schedules should account for the whole range of plausible demand and let the operator decide how much reliability to pay for.} Chance-constrained UC has a long history in power systems 
\cite{ozturk2004solution}. In district heating, uncertainty-aware UC is more recent, with work on robust CHP commitment and dispatch \cite{zugno2016commitment} and chance-constrained co-optimization \cite{polisetty2025chance}. We build on this line of work by
coupling Bayesian probabilistic demand forecasting with risk-aware UC and comparing two risk
paradigms head-to-head on the real Berlin district heating network.

\rev{These paradigms, chance-constrained programming (CC) and conditional
value-at-risk (CVaR) optimization, capture the two basic notions of shortage
risk, frequency and severity: CC bounds how often demand is unmet, CVaR how
severe the worst shortfalls are. While chance-constrained problems are hard in
general, we impose individual chance constraints per time step, which on
equiprobable scenarios reduce to linear constraints with demand fixed to an
empirical quantile, so the CC model keeps the size of the deterministic MIP,
whereas the CVaR model grows with the number of scenarios.}

Our contributions are: (i) the coupling of Bayesian probabilistic demand forecasts, providing full predictive distributions; (ii) the adaptation of two risk paradigms (CC and CVaR) to our mixed-integer UC model for district heating, enabling a controlled comparison; (iii) validation on real operational data from the Berlin network and on generated benchmark instances of varying size, extending our open-source instance generator \cite{generator}; and (iv) an analysis of the price of reliability.

\section{Probabilistic Demand Scenarios}

\rev{For every demand node, a point forecast is given. For the Berlin network, it
is provided by BEW Berliner Energie und W\"arme GmbH; for the synthetic
instances, it is generated as an autoregressive trajectory with realistic
autocorrelation~\cite{buchholz2026benchmarking}. To capture the uncertainty around these forecasts, we fit a Bayesian time-series
model, implemented in PyMC \cite{abril2023pymc}, to each of them and sample
demand scenarios, i.e., plausible trajectories centered on the forecast, from
its posterior-predictive distribution. A Bayesian model is a natural choice
here: it yields a full predictive distribution by construction, propagates
both parameter and observation uncertainty into the scenarios, and remains
interpretable and cheap to fit per node. Demand is modeled in log space as the
sum of an intercept, a linear trend, and Fourier terms for the dominant annual
winter--summer cycle and the weaker daily and weekly cycles of consumer
behavior. A LogNormal likelihood ensures positivity and captures the
right-skewed demand peaks most critical for security of supply. The posterior
is sampled with the No-U-Turn Sampler, a Markov chain Monte Carlo method. Each
scenario is a full trajectory and thus preserves the temporal correlation
between consecutive time steps that matters for minimum up/down times,
ramping, and storage.}

\section{Base model}
\rev{Intuitively, the model decides for each time step which conversion units
are online and how much they produce, how storages are operated, and how much
is traded with external markets, such that demand is covered at minimum cost
within the technical limits of all units.} Formally, the district heating network is a
directed graph on which we pose a mixed-integer unit-commitment formulation.
Nodes correspond to conversion units $i \in I$ (boilers, CHP units, heat
pumps), storage units $k \in K$, transport and emission nodes, and demand
nodes, connected by flows of resources $r \in R$ (heat, fuel, electricity)
over the discretized time horizon $t \in T$. \rev{Commitment is captured by
the binary variables $z_{it}$, with $z_{it} = 1$ if unit $i$ is online in time
step $t$, and $s_{it}$, with $s_{it} = 1$ if unit $i$ is started up in time
step $t$. The continuous flow variables $x$ comprise, among others, the input
and output flows $x^{r}_{t i^{\mathrm{in}}}$ and $x^{r}_{t i^{\mathrm{out}}}$
of resource $r$ at unit $i$ and, analogously, the charging and discharging
flows $x^{r}_{t k^{\mathrm{in}}}$ and $x^{r}_{t k^{\mathrm{out}}}$ of storage
$k$. The quantities purchased from and sold to external markets are $p^{r}_{t}$
and $e^{r}_{t}$. The demand balance requires the demand $d^{r}_{t}$ to be
covered in every time step:}
\begin{equation}
    \sum_{i \in I} x^{r}_{t i^{\mathrm{out}}} + \sum_{k \in K} x^{r}_{t k^{\mathrm{out}}}
      + p^{r}_{t}
    = d^{r}_{t}
      + \sum_{i \in I} x^{r}_{t i^{\mathrm{in}}} + \sum_{k \in K} x^{r}_{t k^{\mathrm{in}}}
      + e^{r}_{t} \qquad \forall r \in R,\, t \in T.
    \label{eq:balance}
\end{equation}
The remaining constraints (conversion, minimum up/down times, ramping, storage
dynamics, and capacity bounds) are detailed in \cite{riedmuller2025enhancing}.
\rev{The objective minimizes the total operating cost $C(x,p,e,z,s)$ over the
horizon, consisting of market purchases net of sales, start-up and running
costs of conversion units, transport costs, and the cost of emission
certificates.}

\section{Risk-aware formulations}

Let $\Omega$ be the set of (equiprobable) scenarios and $d_{\omega t}^r$ be the heat demand
for scenario $\omega \in \Omega$ at time step $t \in T$ for resource $r \in
R$. \rev{Both formulations replace the single demand $d^r_t$ of the base model by this scenario set and differ only in how they treat demand that is not covered.}

\subsection{Chance-constrained programming}
Individual chance-constrained programming requires that demand is met with at least a given
probability $\beta$ at each time step for each resource. Since demand enters as a requirement to
be covered, this is equivalent to sizing supply to the upper $\beta$-quantile
of demand: on the finite, equiprobable scenario
set $\Omega$, this quantile $d^{r}_{\beta t}$ is the $\lceil \beta |\Omega| \rceil$-th order
statistic \rev{(i.e., the $\lceil \beta |\Omega| \rceil$-th smallest value)} of $\{d^{r}_{\omega t}\}_{\omega \in \Omega}$, and $\beta = 1$ recovers the sample worst case. The demand balance \eqref{eq:balance} of the base model is stated with
$d^{r}_{t}$ replaced by the quantile $d^{r}_{\beta t}$. All remaining
base-model constraints stay unchanged, so the resulting model has the same
size and complexity as the deterministic MIP.

\subsection{CVaR optimization}

\rev{Whereas CC bounds how often demand is unmet, CVaR measures how severe the
shortfalls are: for a confidence level $\alpha \in (0,1)$ and a loss random
variable~$L$, $\mathrm{CVaR}_\alpha(L)$ is the mean loss over the worst
$(1-\alpha)$ fraction of outcomes. Rockafellar and
Uryasev~\cite{rockafellar2000optimization} showed that it can be minimized
via an auxiliary threshold variable, which yields a linear formulation for
scenario-based losses.} \rev{Let $u^{r}_{\omega t} \geq 0$ denote the unmet demand of resource $r$ at time $t$ in scenario $\omega$.} We define the loss of scenario $\omega$ as the penalized
unmet demand $L_\omega \;=\; \rho \sum_{t \in T} \sum_{r \in R} u^{r}_{\omega t}$
with unit penalty $\rho$. Introducing the \rev{auxiliary} VaR variable $\zeta \in \mathbb{R}$\rev{, whose optimal value is the loss threshold exceeded with probability $1-\alpha$,} and the scenario excess $c_{\omega} \geq 0$ linearizing $(L_{\omega} - \zeta)^+$, the objective augments the operating cost $C(x,p,e,z,s)$  with the CVaR term, weighted by the risk-aversion parameter $\lambda \geq 0$: 
\begin{equation}
    \min \; C(x,p,e,z,s) + \lambda \Bigl( \zeta + \frac{1}{(1-\alpha)\,|\Omega|}
    \sum_{\omega \in \Omega} c_\omega \Bigr)
    \label{eq:obj}
\end{equation}
subject to the base constraints with the demand balance constraint replaced by:
\begin{align}
    u^{r}_{\omega t} &\geq
      d^{r}_{\omega t}
      + \sum_{i \in I} x^{r}_{t i^{\mathrm{in}}}
      + \sum_{k \in K} x^{r}_{t k^{\mathrm{in}}}
      + e^{r}_{t}
      - \sum_{i \in I} x^{r}_{t i^{\mathrm{out}}}
      - \sum_{k \in K} &&x^{r}_{t k^{\mathrm{out}}}
      - p^{r}_{t}
      \notag \\
    & && \forall r \in R,\ t \in T,\ \omega \in \Omega, \notag \\
    c_\omega &\geq \rho \sum_{t \in T} \sum_{r \in R} u^{r}_{\omega t} - \zeta
      && \forall \omega \in \Omega.
\end{align}

Setting $\lambda = 0$ recovers the penalty-free model. The model is single-stage, with unmet
demand as the only recourse. \rev{In contrast to CC, it contains shortage variables and constraints for every scenario, so its size grows linearly with $|\Omega|$.}

\section{Computational study}
Both methods are implemented in Rust; the resulting MIPs are solved
using Gurobi~13.0 \cite{gurobi} on an Intel Xeon Gold 6338
(2.00\,GHz, 32~threads) with a relative MIP gap of $10^{-3}$ and a
time limit of 24\,h per run. The real-world case study is VG2, the
East-Berlin subgrid of the Berlin district heating network, provided by BEW Berliner Energie und W\"arme GmbH, with the
demand model fitted to six months of forecasting data at eight-hour
resolution, chosen for runtime reasons. The synthetic test bed comprises four instances at
four-hour resolution, each differing from its predecessor in a single
aspect: SYN1 covers a one-year horizon without storage, SYN2 extends the
horizon to two years, SYN3 additionally admits storage, and SYN4 is
posed on a smaller network. All instances use 1\;000 demand scenarios.
The CVaR formulations range from 1.7\,M constraints and 1.7\,M
variables (VG2) to 14.5\,M constraints and 14.2\,M
variables (SYN2). We fix $\alpha = 0.95$, set the
penalty to $\rho = 1.5\,c_{\max}$ with $c_{\max}$ the worst-case
marginal heat cost over all converters and time steps, so that the penalty strictly exceeds generation cost, and sweep
$\beta \in \{0.80, 0.81, \dots, 1.00\}$ and
$\lambda \in \{0, 0.5, 1, 2, \dots, 10\}$.

To compare solutions across scenarios, we define the expected energy not served per year as
\begin{equation}
\text{EENS} \;=\; \frac{1}{|\Omega|} \frac{1}{N}\sum_{\omega\in\Omega}\sum_{t\in T} \sum_{r \in R} \Delta_t
\max\left\{0,\; d^{r}_{\omega t} - q^r_t\right\},
\label{eq:meanshortfall}
\end{equation}
with $q^r_t$ being the realized supply to demand nodes, $\Delta_t$ the temporal resolution and $N$ the number of years. All EENS values are evaluated out of sample: the dispatch obtained on the
optimization scenarios is fixed and evaluated on 1000 newly drawn
scenarios~\cite{mak1999monte}.



\subsection{Chance-constrained programming}
Table~\ref{tab:cc} shows the effect of the reliability level $\beta$
on operational cost and EENS. On the generated instances, cost rises monotonically in $\beta$ while the
out-of-sample EENS decreases. The differences between the instances follow their construction: adding storage (SYN3 vs. SYN2) leaves solutions similar but raises converter uptime and run times. The Berlin instance VG2 behaves monotonically up to $\beta = 0.97$; beyond
that, network capacity no longer covers the quantile demand and the model becomes infeasible. Computational effort for VG2 ranges from 4.9 to 142 minutes, making CC well suited for cheaply tracing
the cost--reliability frontier. 

\begin{table}[!t]
\setlength{\tabcolsep}{3.5pt} 
\caption{Chance-constrained results for selected reliability levels $\beta$.
Cost is reported in M€; EENS (MWh/year) is the out-of-sample shortage.}
\label{tab:cc}
\centering
\begin{tabular}{l rr rr rr rr}
\toprule
 & \multicolumn{2}{c}{$\beta=0.80$} & \multicolumn{2}{c}{$\beta=0.90$}
 & \multicolumn{2}{c}{$\beta=0.95$} & \multicolumn{2}{c}{$\beta=1.00$}\\
Instance & Cost & EENS & Cost  & EENS
         & Cost  & EENS & Cost & EENS\\
\midrule
    VG2  & 531.7   & 178{,}886 & 983.2   & 93{,}201 & 1{,}832.7 & 67{,}587 & \multicolumn{2}{c}{infeasible} \\                                                                                                                        
  SYN1 & 146.1   & 207{,}589 & 169.5   & 97{,}382 & 192.9     & 47{,}057 & 656.5     & 893 \\                                                                                                                                       
  SYN2 & 821.0   & 201{,}159 & 924.5   & 94{,}105 & 1{,}022.8 & 45{,}455 & 1{,}761.0 & 806 \\                                                                                                                                       
  SYN3 & 730.7   & 192{,}504 & 826.6   & 89{,}885 & 918.3     & 43{,}323 & 1{,}499.9 & 763 \\                                                                                                                                       
  SYN4 & 479.0   & 200{,}949 & 553.3   & 93{,}958 & 629.5     & 45{,}353 & 3{,}321.6 & 813 \\
  
\bottomrule
\end{tabular}
\end{table}

\subsection{CVaR optimization}
Table~\ref{tab:cvar} reports some of the CVaR results. On VG2, $\lambda = 1$
leaves the penalty-free schedule essentially unchanged, from
$\lambda = 3$ the CVaR term reshapes commitment and dispatch, and at
$\lambda = 10$ cost more than triples while the mean shortage falls by
$90\,\%$. SYN1 responds threshold-like: at $\lambda = 0$ shedding
demand is profitable, while $\lambda = 1$ already cuts the shortage by
an order of magnitude. Computational effort scales with $\lambda$, at $\lambda = 10$ instance SYN1 did not converge within the 24 hours.

\begin{table}[!t]
\centering
\setlength{\tabcolsep}{6.5pt} 
\caption{CVaR optimization results of VG2 and SYN1 for selected risk-aversion weights $\lambda$. Operational cost and CVaR are reported in M\texteuro; EENS (MWh/year) is the out-of-sample shortage.}
\begin{tabular}{lrrrrcr}
\toprule
Instance & $\lambda$ & Cost & CVaR & Time (s) & MIP gap & EENS \\
\midrule
VG2 & 0.0  & 92.4  & 53.8  & 7      & optimal & 5,524,384 \\
VG2 & 1.0  & 92.4  & 53.6  & 245    & optimal & 5,516,393 \\
VG2 & 3.0  & 106.4 & 47.8  & 86,400 & 0.0046  & 4,906,522 \\
VG2 & 10.0 & 303.1 & 6.3   & 86,400 & 0.0124  & 565,886 \\
\hline
SYN1 & 0.0  & -6.8   & 371.4 & 222    & optimal & 5,074,590 \\
SYN1 & 1.0  & 118.8  & 33.1  & 86,400 & 0.0019  & 417,928 \\
SYN1 & 3.0  & 153.6  & 11.3  & 86,400 & 0.0012  & 136,279 \\
\bottomrule
\end{tabular}
\label{tab:cvar}
\end{table}

\begin{figure}[t]
  \centering
  \includegraphics[width=0.54\linewidth]{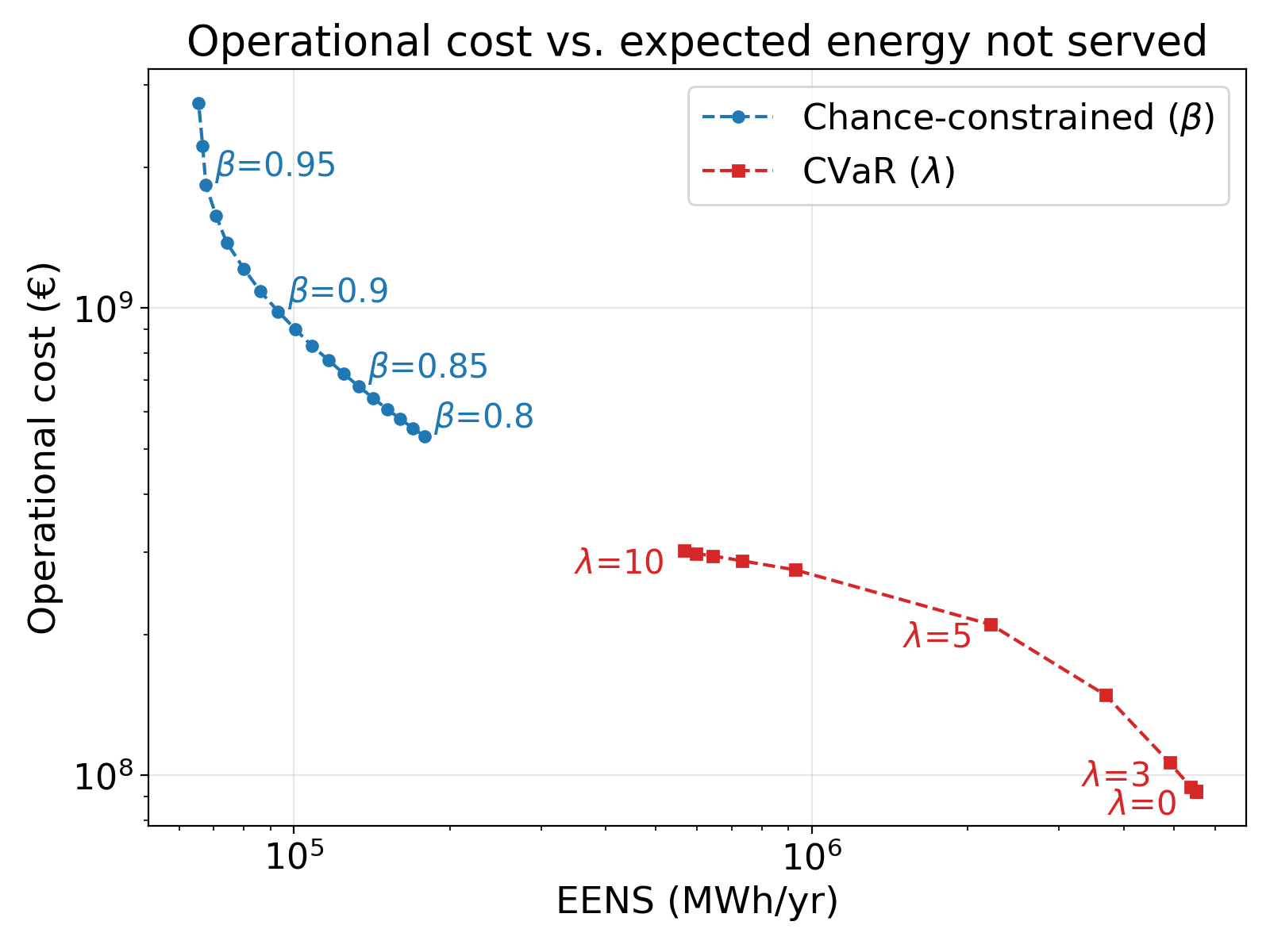}
  \caption{Cost--reliability of both methods on VG2
  (log--log): out-of-sample EENS (MWh/year)
  versus operational cost (€) (excluding the CVaR penalty term),
  for the chance-constrained sweep $\beta \in \{0.80,\dots,0.97\}$
  and the CVaR sweep $\lambda \in \{0,\dots,10\}$.}
  \label{fig:frontier}
\end{figure}

\subsection{Comparison}
Figure~\ref{fig:frontier} places both methods in the same
cost--reliability plane, out of sample. CC schedules reduce out-of-sample shortage by roughly a factor of ten compared to the most risk-averse CVaR run, at several times the cost, while the CVaR
sweep covers the inexpensive, risk-tolerant regime the $\beta$-sweep
does not. Neither dominates: the two methods control different risk objects, shortfall probability per time step (CC) versus tail severity of shortfall (CVaR), and thus trace complementary segments of the frontier.

\section{Conclusion}
\rev{Our results show that reliability has a steep and increasing price: on VG2,
raising $\beta$ from $0.80$ to $0.95$ more than triples the operating cost
while the expected energy not served falls by about $60\,\%$, and on the
synthetic instances, covering the sample worst case costs $1.6$ to $5.3$ times
as much as $\beta = 0.95$. Making this trade-off explicit lets operators choose
a reliability level deliberately instead of implicitly through safety margins
on a point forecast. CC and CVaR do not compete but cover complementary parts
of this trade-off. In practice, this suggests a clear division of labor. CC fits when
reliability is specified as a target: $\beta$ directly states the accepted
shortfall probability per time step, and a sweep, fast enough for routine planning, shows
what each level costs and which levels the network can reach at all. On VG2,
infeasibility beyond $\beta = 0.97$ signals that, under the generated demand
scenarios, higher reliability requires investment in capacity rather than a
better schedule. CVaR fits when the size
of a shortfall matters more than its occurrence, or when no reliability target
is attainable: since unmet demand is penalized rather than forbidden, CVaR
always returns a schedule, and $\lambda$ expresses how much the operator is
willing to pay to avoid severe shortfalls, albeit at considerably higher
computational cost. The two can thus be combined: a CC sweep first maps the
achievable reliability and its price, and CVaR then shapes schedules in the
risk-tolerant regime or where the reliability target cannot be met.}
Future work includes \rev{joint chance constraints over the horizon,}
uncertainty in fuel, electricity, and CO$_2$ prices, a multi-objective
formulation minimizing emissions, and a two-stage model fixing commitment in
the first stage.

\section*{Acknowledgement}
The work for this article has been conducted in the Research Campus MODAL funded by the German Federal Ministry of Research, Technology and Space (BMFTR) (fund numbers 05M14ZAM, 05M20ZBM, 05M2025).\\
Code used within this manuscript has been developed with the assistance of artificial intelligence (AI) using Claude Code with Claude Opus 4.8 developed by Anthropic (Claude Opus 4.8, June–July 2026). The generated content was carefully reviewed, verified, and endorsed by the authors to ensure the accuracy and integrity of the final manuscript.

\bibliographystyle{unsrt}
\bibliography{GOR/refs}

@article{ozturk2004solution,
  title={A solution to the stochastic unit commitment problem using chance constrained programming},
  author={Ozturk, U. A. and Mazumdar, M. and Norman, B. A.},
  journal={IEEE Transactions on Power Systems},
  volume={19},
  number={3},
  pages={1589--1598},
  year={2004},
  publisher={IEEE}
}

@article{zugno2016commitment,
  title={Commitment and dispatch of heat and power units via affinely adjustable robust optimization},
  author={Zugno, M. and Morales, J. M. and Madsen, H.},
  journal={Computers \& Operations Research},
  volume={75},
  pages={191--201},
  year={2016},
  publisher={Elsevier}
}

@article{polisetty2025chance,
  title={Chance constrained co-optimization of integrated electrical and district heating networks},
  author={Polisetty, S. P. and Nazir, F. U. and Pal, B. C.},
  journal={IEEE Transactions on Power Systems},
  volume={40},
  number={4},
  pages={3401--3412},
  year={2025},
  publisher={IEEE}
}

@article{rockafellar2000optimization,
  title={Optimization of conditional value-at-risk},
  author={Rockafellar, R. T. and Uryasev, S.},
  journal={Journal of risk},
  volume={2},
  pages={21--42},
  year={2000}
}

@misc{gurobi,
  author = {{Gurobi Optimization, LLC}},
  title  = {Gurobi Optimizer Reference Manual},
  year   = {2025},
  note   = {Available at: \url{https://www.gurobi.com}}
}

@misc{generator,
  author       = {Buchholz, A. and Riedm{\"u}ller, S. and Passage, M. and Zittel, J.},
  title        = {Multi Energy Instance Generator},
  howpublished = {https://git.zib.de/abuchhol/multi-energy-instance-generator},
  note         = {Accessed: 31.07.2026}
}

@article{abril2023pymc,
  title={PyMC: a modern, and comprehensive probabilistic programming framework in Python},
  author={Abril-Pla, O. and Andreani, V. and Carroll, C. and Dong, L. and Fonnesbeck, C. J. and Kochurov, M. and Kumar, R. and Lao, J. and Luhmann, C. C. and Martin, O. A. and others},
  journal={PeerJ Computer Science},
  volume={9},
  pages={e1516},
  year={2023},
  publisher={PeerJ Inc.}
}

@inproceedings{riedmuller2025enhancing,
  author    = {Riedm{\"u}ller, S. and Buchholz, A. and Zittel, J.},
  title     = {Enhancing Multi-Energy Modeling: The Role of Mixed-Integer Optimization Decisions},
  booktitle = {Proceedings of the 38th International Conference on Efficiency, Cost, Optimization, Simulation and Environmental Impact of Energy Systems (ECOS 2025)},
  year      = {2025}
}

@inproceedings{buchholz2026benchmarking,
  author        = {Buchholz, A. and Riedm{\"u}ller, S. and Passage, M. and Zittel, J.},
  title         = {Benchmarking Realistic Synthetic Instances Against a Large-Scale District Heating Network: A Multi-Objective Optimization Study for {B}erlin},
  eprint        = {2606.02195},
  archivePrefix = {arXiv},
  primaryClass  = {math.OC},
  howpublished  = {https://arxiv.org/abs/2606.02195},
  booktitle = {Proceedings of the 39th International Conference on Efficiency, Cost,
               Optimization, Simulation and Environmental Impact of Energy Systems (ECOS 2026)},
  year      = {2026},
  note      = {To appear. arXiv:2606.02195},
}

@article{mak1999monte,
  title={Monte Carlo bounding techniques for determining solution quality in stochastic programs},
  author={Mak, Wai-Kei and Morton, David P and Wood, R Kevin},
  journal={Operations research letters},
  volume={24},
  number={1-2},
  pages={47--56},
  year={1999},
  publisher={Elsevier}
}

\end{document}